\documentclass[11pt]{amsart}
\usepackage{amsmath,amssymb,amsfonts,amscd,graphicx}
\usepackage{color}
\usepackage{latexsym}
\usepackage{amsmath}
\usepackage[all]{xy}
\usepackage{hyperref}
\hypersetup{
    colorlinks = true,
    linkcolor = {black},
   citecolor  = {black},
}
\newcommand{\co}{\colon\thinspace}

\numberwithin{equation}{section}

\newtheorem{cor}{Corollary}[section]

\newtheorem{theoremsection}[cor]{Theorem}

\theoremstyle{definition}
\newtheorem{defi}[cor]{Definition}

\theoremstyle{remark}

\begin{document}

\title[TQFT and polynomial identities for graphs on surfaces]{Topological quantum field theory and \\polynomial identities for graphs on the torus}

\author{Paul Fendley and Vyacheslav Krushkal}

\address{Paul Fendley{\hfil\break} All Souls College and Rudolf Peierls Centre for Theoretical Physics, Parks Rd, Oxford OX1 3PU, UK}
\email{paul.fendley@physics.ox.ac.uk}

\address{Vyacheslav Krushkal{\hfil\break} Department of Mathematics, University of Virginia,
Charlottesville, VA 22904-4137 USA}
\email{krushkal\char 64 virginia.edu}

 \begin{abstract} 
We establish a relation between the trace evaluation in ${\rm SO}(3)$ topological quantum field theory and evaluations of a topological Tutte polynomial. As an application, a generalization  of the Tutte golden identity is proved  for graphs on the torus.
\end{abstract}

\maketitle

\section{Introduction}

The Witten-Reshetikhin-Turaev topological quantum field theory (TQFT) associates  invariants to ribbon graphs in $3$-manifolds. A part of this theory is an invariant of graphs on surfaces: given a graph $G\subset \Sigma$, the {\em trace evaluation} is the invariant associated to the embedding $G\subset {\Sigma}\times \{*\}\subset  {\Sigma}\times S^1$. We study the trace evaluation for ${\rm SO}(3)$ TQFTs.

For {\em planar} graphs, the ${\rm SO}(3)$ quantum evaluation is known to equal the flow polynomial $F^{}_G$, or equivalently the chromatic polynomial ${\chi}^{}_{G^*}$ of the dual graph. In \cite{FK} we showed that this quantum-topological approach gives a conceptual framework for analyzing the relations satisfied by the chromatic and flow polynomials of planar graphs. In particular, we gave a proof in this setting of the Tutte golden identity \cite{T2}: given a planar triangulation $T$,
\begin{equation} \label{golden identity eq}
  {\chi}^{}_T({\phi}+2)=({\phi}+2)\; {\phi}^{3\,V(T)-10}\,
  ({\chi}^{}_T({\phi}+1))^2,
\end{equation}
where $V(T)$ is the number of vertices of the triangulation, and ${\phi}$ denotes the golden ratio, ${\phi}=\frac{1+\sqrt5}{2}$. The dual formulation in terms of the flow polynomial states that for a planar cubic graph $G$, $F^{}_G({\phi}+2) \, = \,  {\phi}^E\, (F^{}_G({\phi}+1))^2.$ In \cite{FK} we also showed that the golden identity may be thought of as a consequence of level-rank duality between the ${\rm SO}(3)_4$ and the ${\rm SO}(4)_3$ TQFTs, and the isomorphism ${\mathfrak{so}}(4)\cong {\mathfrak{so}}(3)\times {\mathfrak{so}}(3)$. (Consequences of ${\rm SO}$ level-rank duality for link polynomials have also been studied in \cite{MPS}.)

The main purpose of this paper is to formulate an extension of the results on TQFT and polynomial invariants to graphs on the torus; 
in particular we prove a generalization of the golden identity. Some of the motivation for this work has its origins in lattice models in statistical mechanics, where it has long been known (see e.g.\ \cite{DSZ}, \cite{Pasquier}) that when deriving identities for partition functions  {on the torus}, one must typically sum over ``twisted'' sectors. Twisted sectors are more complicated analogs of the spin structures  (cf. \cite{CR}) familiar in field theories involving fermions. Such sectors are described naturally in TQFT, as for example can be seen in the study of lattice topological defects \cite{AMF}.  We show how the golden identity is generalised to the torus precisely by considering such sums over analogous sectors. In a forthcoming paper \cite{FK3} we will elaborate further on the connections to statistical mechanics, in particular on the relation with the Pasquier height model  \cite{Pasquier}. 

The chromatic and flow polynomials are $1$-variable specializations of the Tutte polynomial, known in statistical mechanics as the partition function of the Potts model. Relations between the ${\rm SO}(3)$ quantum evaluation of planar graphs, the chromatic and flow polynomial, and the Potts model are discussed in detail in \cite{FK2}. 
From the TQFT perspective, the case of graphs embedded in the plane (or equivalently in the $2$-sphere $S^2$)  is very special in that the
TQFT vector space associated to $S^2$ is ${\mathbb C}$. For surfaces $\Sigma$ of higher genus they are vector spaces (of dimension given by the Verlinde formula) which are part of the rich structure given by the $(2+1)$-dimensional TQFT. Multi-curves, and more generally graphs embedded in $\Sigma$, act as ``curve operators'' on the TQFT vector space, and our goal is to analyze the trace of these operators.

For the ${\rm SO}(3)$ TQFT, this invariant of graphs on surfaces satisfies the contraction-deletion rule, familiar from the study of the Tutte polynomial and the Potts model. A non-trivial feature on surfaces of higher genus is that the ``loop value'' depends on whether the loop is trivial (bounds a disk) in the surface, or whether it wraps non-trivially around the surface; in fact in the latter case the invariant is not multiplicative under adding/removing a loop. This behavior is familiar in  generalizations of the Tutte polynomial which encode the topological information of the graph  embedding in a surface. The study of such ``topological'' graph polynomials was pioneered by Bollob\'{a}s-Riordan \cite{BR}; a more general version was introduced in \cite{Kr}. To express the ${\rm SO}(3)$ trace evaluation for graphs on the torus we need a further extension of the polynomial, defined in section \ref{graph poly section}. In Theorem \ref{evaluations are equal} we show that the quantum invariant equals a sum of values of the polynomial, where the sum is parametrized by labels corresponding to the TQFT level. Here individual summands (evaluations of the topological Tutte polynomial) correspond to TQFT sectors. 

The identity (\ref{golden identity eq}) for the chromatic polynomial and its analogue for the flow polynomial in general do {\em not} hold for non-planar graphs. In fact, it is conjectured \cite{AK} that the golden identity for the flow polynomial characterizes planarity of cubic graphs.  Our result (Theorem \ref{golden torus theorem}) involves invariants of graphs derived from TQFT, or equivalently evaluations of  graph polynomials, and it recovers the original Tutte's identity for planar graphs.

A different ``quantum'' version of the golden identity, for the Yamada polynomial of ribbon cubic graphs in ${\mathbb R}^3$, was established in \cite{AK2}. It is an extension of the Tutte identity from the flow polynomial of planar graphs to the Yamada polynomial of spatial graphs, which may be thought of as elements of the skein module of $S^2$ (isomorphic to ${\mathbb C}$). Our result manifestly involves graphs on the torus and the elements they represent in the associated TQFT vector spaces. We expect that an extension of our results holds on surfaces of genus $>1$ as well, although computational details on surfaces of higher genus are substantially more involved.

To date, the study of topological Tutte polynomials of graphs on surfaces has been carried out primarily in the context of topological combinatorics. While there are applications  to quantum invariants of links (cf. \cite{DFKLS, MM}) and to noncommutative quantum field theory \cite{KRTW}, to our knowledge the results of this paper provide the first direct relation between TQFT and evaluations of graph polynomials on surfaces.

After briefly reviewing background information on TQFT in section \ref{TQFT section}, we define the relevant topological Tutte polynomial for graphs on the torus in section \ref{graph poly section}. The relation between TQFT trace and evaluation of the graph polynomial is stated and proved in section \ref{relation section}. The generalization of the golden identity for graphs on the torus is established in section \ref{golden torus section}.

\section{Background on TQFT} \label{TQFT section}
We refer the reader to \cite{FK} for a summary of the relevant facts about the Temperley-Lieb algebras ${\rm TL}_n$ and the Jones-Wenzl projectors  $p_n$ \cite{Jo, We}, and to \cite{KL} for a more general introduction to the calculus of quantum spin networks.

To fix the notation, recall the definition of quantum integers $[n]$,  and the evaluation of the $n$-colored unknot, ${\Delta}_n$:
\[
[n]=\frac{A^{2n}-A^{-2n}}{A^2-A^{-2}}\ , \qquad  {\Delta}_n=[n+1]\ .
\]

We will interchangeably use the parameters $q, A$, as well as the loop value $d$ and a graph parameter $Q$, 
related as follows:
\begin{align} 
\label{notation1} 
q=A^4, \quad d=A^2+A^{-2}, \quad Q=d^2\ .
\end{align}

A basic calculation for the colored Hopf link gives 
\begin{equation}  \label{LinkingLoop fig}
\vcenter{\hbox{\includegraphics[height=2cm]{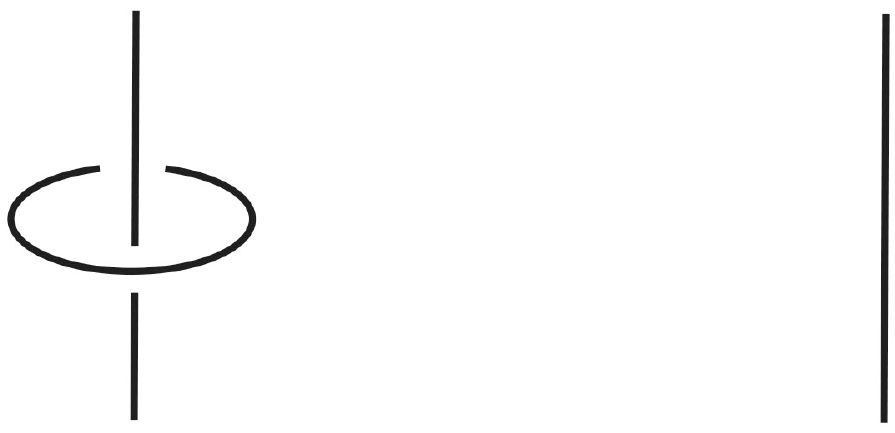}
\put(-131,25){$m$}
\put(-96,50){$n$}
{ \put(-78,22){$= \, \frac{[(m+1)(n+1)]}{[n+1]}$} }
\put(2,50){$n$}
}}
\end{equation}

We use the construction of ${\rm SU}(2)$ Witten-Reshetikhin-Turaev TQFTs, given in \cite{BHMV}. Given a closed orientable surface $\Sigma$, the ${\rm SU}(2)$, level $r-2$ TQFT vector space will be denoted $V_r({\Sigma})$, where $A=e^{2{\pi} i/4r}$. (Note that in the TQFT literature the notation $V_{2r}$ is sometimes used instead. Also note that in the physics literature, labels are often divided by 2 and called ``spin'', so that odd and even labels correspond to half-integer and integer spins respectively.)

The main focus of this paper is on the torus case, ${\Sigma}={\mathbb T}$, and in this case the notation $V_r:=V_r({\mathbb T})$ will be used throughout the paper. Consider ${\mathbb T}$ as the boundary of the solid torus $H$. $V_r$ has a basis $\{ e_0, \ldots, e_{r-2}\}$, where $e_j$ corresponds to the core curve of $H$, labeled by the $j$-th projector $p_j$. This basis will be used in the evaluation of the trace in Section \ref{trace subsection}.

The discussion in the rest of this section applies to surfaces $\Sigma$ of any genus.
A curve $\gamma$ in $\Sigma$ acts as a linear operator on $V_r({\Sigma})$, so associated to $\gamma$ is an element of $V^*_r({\Sigma})\otimes V_r({\Sigma})$.
Given a graph $G\subset \Sigma$, we consider it as an ${\rm SU}(2)$ quantum spin network in $\Sigma$ by turning each edge into a ``double line''. Namely, as in \cite{FK} we label edges by the second Jones-Wenzl projector, up to an overall normalization. 
Concretely,  each edge $e$ of $G$ is replaced with a linear combination of curves as indicated in (\ref{phi}), and the curves are connected without crossings on the surface near each vertex. 
\begin{equation}  \label{phi}
\vcenter{\hbox{
\includegraphics[width=12cm]{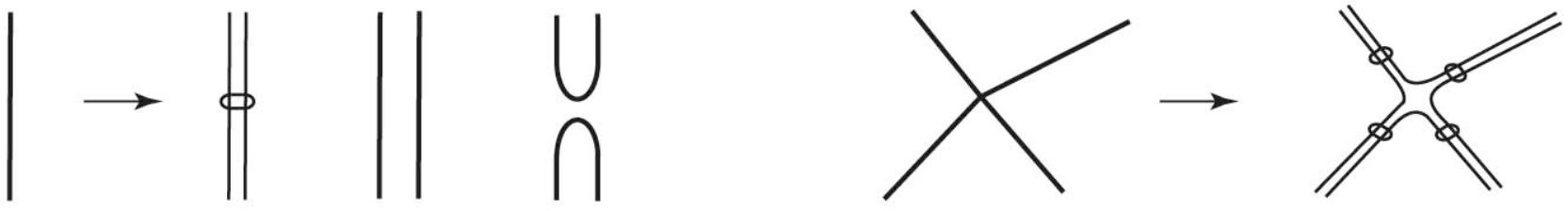}
\put(-350,5){$e$}
    \put(-322,29){${\Phi}$}
    \put(-277,21){$=$}
    \put(-244,21){$-\frac{1}{d}$}
    \put(-85,29){${\Phi}$}
    \put(-62,21){$d\;\cdot$}
}}
\end{equation}
In this map, a factor
$d^{(k-2)/2}$ is associated to each $k$-valent vertex, so that for example the 
$4$-valent vertex in (\ref{phi}) is multiplied by $d$. The overall factor for a graph $G$ is the product of the factors $d^{(k(V)-2)/2}$ over all vertices $V$ of $G$. Therefore the total exponent equals half the sum of valencies over all vertices, minus the number of vertices, i.e.\ minus the Euler characteristic of $G$. (This count does not involve faces - so this is the Euler characteristic of the graph $G$, and not of the underlying surface $\Sigma$.) Using this map $\Phi$, the graph is mapped to a linear combination of multicurves in the surface $\Sigma$. We thus may consider graphs $G$ on the torus as elements ${\Phi}(G)\in{\rm Hom}(V_r, V_r)$.

\begin{figure}[ht]
\includegraphics[height=2cm]{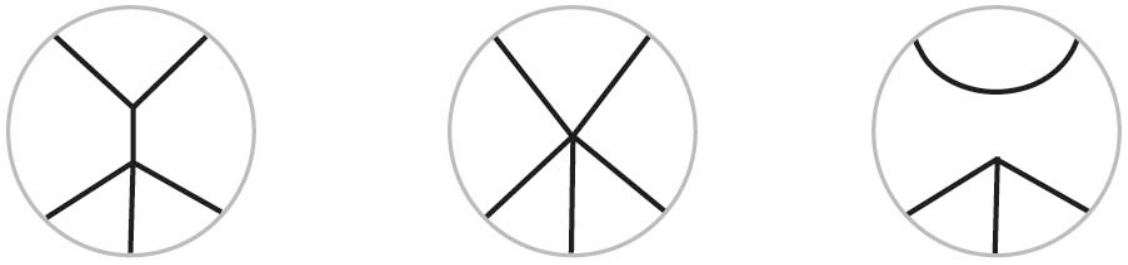}
    \put(-175,28){$=$}
    \put(-80,28){$-$}
{\scriptsize
    \put(-212,26){$e$}
    \put(-245,54){$G$}
    \put(-152,54){$G/e$}
    \put(-7,54){$G\backslash e$}}
\caption{Relation (1)}
\label{fig:rel1}
\end{figure}

\begin{figure}[ht]
\includegraphics[height=2cm]{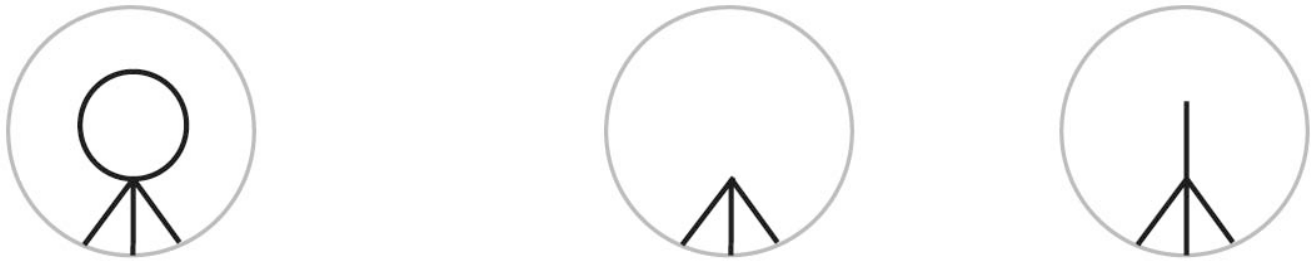}
    \put(3,25){$=\; 0.$}
    \put(-218,25){$=\; (Q-1)\;\cdot$}
    \put(-100,15){$,$}
    \put(-275,33){$e$}
    \put(-39,24){$e$}
\caption{Relations (2), (3)}
\label{fig:rel23}
\end{figure}

Given a graph $G\subset{\mathbb T}$, consider the following local relations (1)--(3), illustrated in figures \ref{fig:rel1}, \ref{fig:rel23}.

(1) If $e$ is an edge of a graph $G$ which is not a loop, then
$G=G/e-G\backslash e$, as illustrated in figure \ref{fig:rel1}.

(2) If $G$ contains an edge $e$ which is a loop, then $G=(Q-1)\;
G\backslash e$, as in figure \ref{fig:rel23}. (In particular,
this relation applies if $e$ is a loop not connected to
the rest of the graph.)

(3) If $G$ contains a $1$-valent vertex 
as in figure \ref{fig:rel23}, then $G=0$.

Replacing the edge labeled $e$ in each figure with the linear combination of curves, defined by the second JW projector, one checks that the relations (1)--(3) hold for ${\Phi}(G)$ in ${\rm Hom}(V_r, V_r)$, where the parameters $A, Q$ are related as in (\ref{notation1}), and the value of $A$ in the definition of $V_r$ is $e^{2{\pi} i/4r}$. 

{\bf Remark.}
In the planar case, the map gives a homomorphism from the {\em chromatic algebra} ${\mathcal
  C}^{Q}_n$, defined in \cite{FK}, to the Temperley-Lieb algebra $TL^d_{2n}$, with the parameters related by $Q=d^2$. This map is used to show \cite[Lemma 2.5]{FK} that for planar graphs $G$, up to a normalization the quantum evaluation equals the flow polynomial of $G$ or equivalently the chromatic polynomial of the dual graph $G^*$. Indeed,  the relations (1)--(3) are sufficient for evaluating any graph {\em in the plane}, and these relations are precisely the defining relation for the chromatic polynomial (of the dual graph). On the torus, or any other surface of higher genus, the ``evaluation'' is not an element of ${\mathbb C}$ but rather an element of the higher-dimensional TQFT vector space ${\rm Hom}(V_r, V_r)$.

A crucial feature underlying the construction of TQFTs is that for each $r$, {\em in addition to} (1)--(3) there is another local relation corresponding to the vanishing of the corresponding Jones-Wenzl projector. For example, consider the case $r=5$, important in the proof of the golden identity below. In this case, graphs G$\subset {\Sigma}$, considered as elements of the vector space ${\rm Hom}(V_{5}({\Sigma}), V_{5}({\Sigma}))$,  satisfy the local relation in  (\ref{fig:graphp4}).
\begin{equation}  \label{fig:graphp4}
\vcenter{\hbox{
\includegraphics[height=1.6cm]{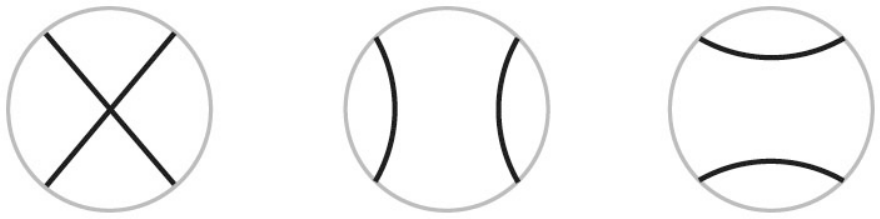}
{\small
    \put(-193,20){$\phi$}
    \put(-130,20){$=$}
    \put(-60,20){$+$}}
{\tiny
    \put(-142,2){$G^{}_X$}
    \put(-72,2){$G^{}_I$}
    \put(-5,2){$G^{}_E$}
    }
}}
\end{equation}
This relation (discovered in the setting of the chromatic polynomial of planar graphs by Tutte \cite{T1}) corresponds to the $4$th JW projector, see \cite[Section 2]{FK} for more details. 

{\bf Remark.} 
In fact, the Turaev-Viro ${\rm SU}(2)$ TQFT associated to a surface $\Sigma$, isomorphic to ${\rm Hom}(V_{r}({\Sigma}), V_{r}({\Sigma}))$, can be defined as the vector space spanned by multi-curves on $\Sigma$, modulo the local relations given by ``$d$-isotopy'' and the vanishing of the JW projector. (See  Theorem 3.14 and section 7.2 in \cite{FNWW}.)
The ${\rm SO}(3)$ theory may be built by considering the ``even labels'' subspace spanned by graphs modulo relations (1)--(3) above, and the JW projector.

\section{Polynomial invariants of graphs on surfaces} \label{graph poly section}

Let $H$ be a graph embedded in the torus ${\mathbb T}$. 
Let ${\rm n}(H)$ denote the nullity of $H$, that  is the rank of the first homology group $H_1(H; {\mathbb Z})$.
The rank $r(H\subset {\mathbb T})$ of the image of the map $i_*\co H_1(H; {\mathbb Z})\longrightarrow H_1({\mathbb T}; {\mathbb Z})\cong {\mathbb Z}^2$, induced by the inclusion $H\subset {\mathbb T}$, is either $0$, $1$ or $2$.
Consider the following homological invariants: 
\begin{equation}
 {\rm s}(H) :=
  \begin{cases}
    1 & \text{if $r(H\subset {\mathbb T}) =2 $ (in other words, $H_1(H)\longrightarrow H_1({\mathbb T})$ is surjective),} \\
    0 & \text{otherwise}.
  \end{cases}
\end{equation}
\begin{equation}
 {\rm s}^{\perp}(H) :=
  \begin{cases}
    1 & \text{if $r(H\subset {\mathbb T}) =0 $ (i.e. $H_1(H)\longrightarrow H_1({\mathbb T})$ is the zero map),} \\
    0 & \text{otherwise.}
  \end{cases}
\end{equation}
Let ${\rm c}(H)$ be the number of rank $0$ connected components of $H$, that is the number of connected components $H^{(i)}$ such that $r(H^{(i)}\subset {\mathbb T})=0$.

In case $r(H\subset {\mathbb T}) =1$, let $\bar {\rm c}(H)$ be the number of ``essential'' components of $H$, i.e. the number of components $H^{(i)}$ such that $H_1(H^{(i)})\longrightarrow H_1({\mathbb T})$ is non-trivial. (Note that for each such component, the image of the map on homology is the same rank $1$ subgroup of $H_1({\mathbb T})$.) If  $r(H\subset {\mathbb T})$ is $0$ or $2$, $\bar {\rm c}(H)$ is defined to be zero.

The most general polynomial, encoding the homological information of the embedding of a graph $G$ in the torus, is defined by the following state sum:
\begin{equation} \label{general poly definition}
\widetilde{P}^{}_G(X,Y,W,A,B)\, :=\, \sum_{H\subset G}  (-1)^{E(G)-E(H)} X^{{\rm c}(H)}\,  Y^{{\rm n}(H)}\,  W^{\bar {\rm c}(H)}\,  A^{{\rm s}(H)}\,  B^{{\rm s}^{\perp}(H)},
\end{equation}
where the summation is taken over all spanning subgraphs of $G$, $E(G)$ denotes the number of edges of $G$, and $(-1)^{E(G)-E(H)}$ provides a convenient normalization. Note that our convention for the sign and the variables differs from the usual convention for the Tutte polynomial.
The usual proof (cf. \cite[Lemma 2.2]{Kr}) shows that this polynomial satisfies the contraction-deletion relation.

{\bf Remark}. The polynomials defined in  \cite{BR, Kr} are  specializations of $\widetilde P$ in the case of graphs embedded in the torus.

To establish a relation with the trace evaluation in TQFT, consider the specialization of $\widetilde{P}$ obtained by setting $X=B=1$:
\begin{equation} \label{poly definition}
P^{}_G(Y,W,A)\, :=\, \sum_{H\subset G} (-1)^{E(G)-E(H)} \,  Y^{{\rm n}(H)}\,  W^{\bar {\rm c}(H)}\,  A^{{\rm s}(H)}.
\end{equation}
 This is a generalization of the flow polynomial, including variables $W$ and $A$ which reflect the topological information of how the graph $G$ wraps around the torus. In particular, if $G$ is homologically trivial on the torus (the rank $r(H\subset {\mathbb T})$ is zero), $P^{}_G(Y,W,A)$ recovers the flow polynomial $F^{}_G(Y)$. We thus name $P_G$ the ``topological flow polynomial''.

\begin{figure}[ht] 
\includegraphics[height=2.5cm]{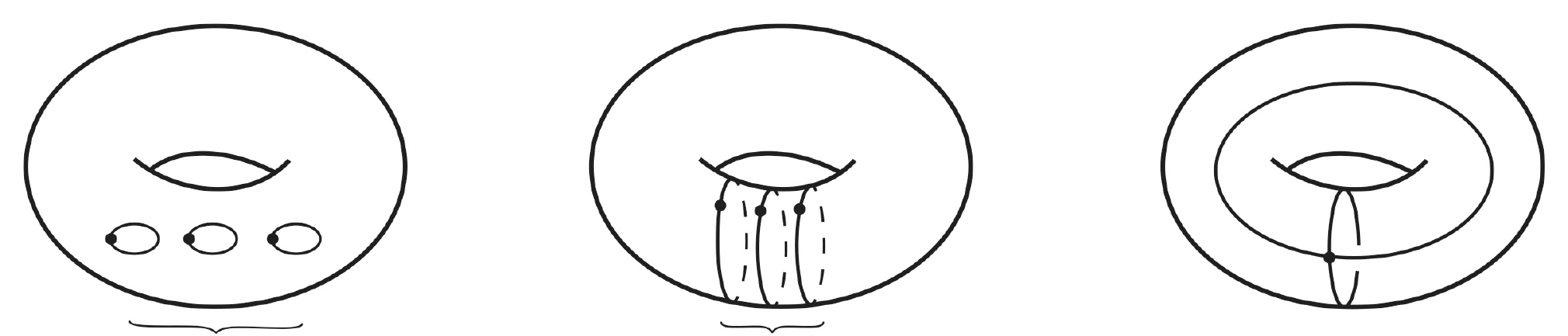}
\put(-290,-14){$k$}
\put(-172,-14){$k$}
\put(-335, -5){(a)}
\put(-215, -5){(b)}
\put(-95,-5){(c)}
\caption{}
\label{trace11 fig}
\end{figure}

Three key examples are illustrated in figure \ref{trace11 fig}. 
For the graph consisting of $k$ disjoint, trivial loops on the torus in figure \ref{trace11 fig}a, the polynomial $P=(Y-1)^k$.
For the graph consisting of $k$ non-trivial loops in figure \ref{trace11 fig}b,  $P=(YW-1)^k$. The calculation is simple; for example the subgraphs $H$ are
$$ 
\vcenter{\hbox{\includegraphics[height=2.5cm]{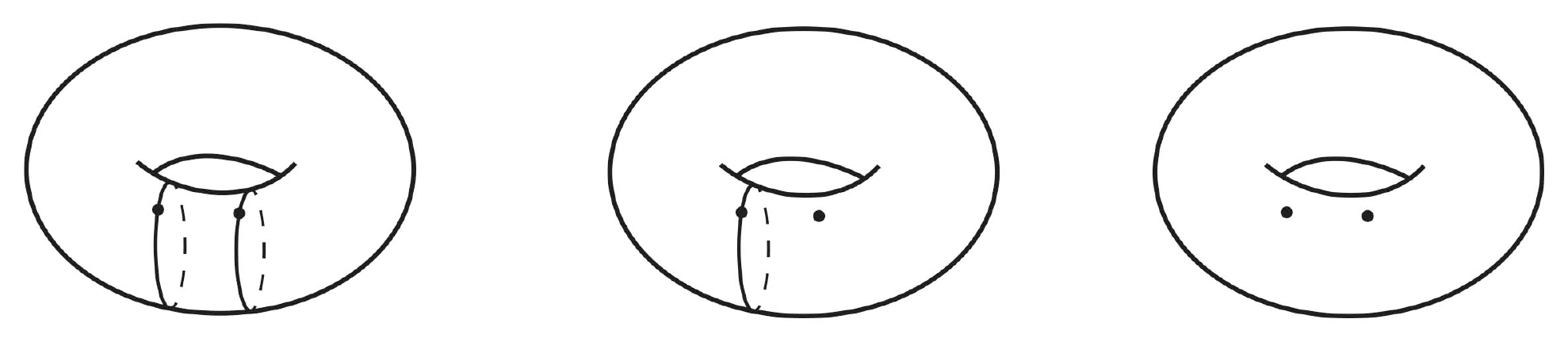}
\put(-230,30){$+\; 2$}
\put(-108,30){$+$}
}}
$$
and so here the polynomial is
\begin{equation}  \label{trace12 fig}=\; Y^2W^2-2YW+1\; =\; (YW-1)^2.
\end{equation}
Finally, the polynomial of the graph  in figure \ref{trace11 fig}c is computed as follows:
$$
\vcenter{\hbox{
\includegraphics[height=2.5cm]{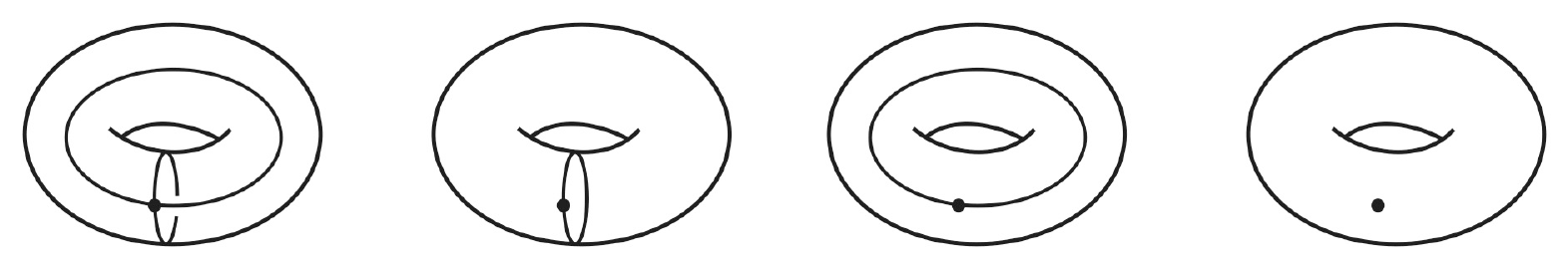}
\put(-321,30){$+$}
\put(-214,30){$+$}
\put(-107,30){$+$}
}}
$$
\begin{equation} \label{poly example}
=\; AY^2- 2YW+1.
\end{equation}

\subsection{Duality} \label{duality section}
The chromatic and flow polynomials ${\chi}, F$ are $1$-variable specializations of the Tutte polynomial, satisfying the relation $F^{}_G(Q)=Q^{-1}{\chi}^{}_{G^*}(Q)$, where $G$ is a planar graph and $G^*$ is its dual. We can extend this duality to the torus
by defining the ``topological chromatic polynomial'. Namely, we specialize $\widetilde{P}_G$ in \eqref{general poly definition} to
\begin{equation} \label{chrom poly definition}
{C}^{}_G(X,U,B)\, :=\, \sum_{H\subset G}  (-1)^{E(G)-E(H)} X^{{\rm c}(H)-{\rm c}(G)}\,   U^{\bar {\rm c}(H)}  B^{{\rm s}^{\perp}(H)}.
\end{equation}
When $r(H\subset {\mathbb T}) =0$, ${C}^{}_G(X,U,B)$ is a renormalized version of the chromatic polynomial ${\chi}^{}_G(X)$. Analogously to the proof of \cite[Theorem 3.1]{Kr}, one shows that for a cellulation $G\subset {\mathbb T}$ (i.e. when each face of the embedding is a $2$-cell),
\begin{equation} \label{duality eq}
{P}^{}_G(Y,W,A)  \, =\,  {C}^{}_{G^*}(Y,YW,AY^2),
\end{equation}
where $G^*\subset {\mathbb T}$ is the dual graph. This topological chromatic polynomial is thus dual to the topological flow polynomial $P_G$ from (\ref{poly definition}).

\section{TQFT trace as an evaluation of  a graph polynomial} \label{relation section}

 \subsection{TQFT trace} \label{trace subsection}

 Given an odd $r$ and a multi-curve $\gamma\subset {\mathbb T}$, the trace of $\gamma$ is defined as $Z_r({\mathbb T}\times S^1, {\gamma})$, the quantum invariant of the banded link $\gamma$ in the $3$-manifold ${\mathbb T}\times S^1$ \cite[1.2]{BHMV}. Concretely, it may be calculated as the trace of the curve operator in ${\rm Hom}( V_r, V_r)$ with respect to the basis discussed in Section \ref{TQFT section}. This basis is given by the core circle of a fixed solid torus $H$, bounded by ${\mathbb T}$, labeled with an integer $0\leq j\leq r-2$. For a graph $G\subset{\mathbb T}$, the trace is defined by mapping $G$ to a linear combination of multicurves using ${\Phi}$ in (\ref{phi}) and then computing the trace of ${\Phi}(G)$.
 
 While one could work with the full ${\rm SU}(2)$ TQFT vector space $V_r$, the invariants of graphs obtained by putting the $2$nd Jones-Wenzl projector on the edges, as in (\ref{phi}), naturally fit in the context of the ${\rm SO}(3)$ theory. This corresponds to taking the subspace $\overline V_r$ of $V_r$, spanned by even labels. For the remainder of the paper, ${\rm tr}_r(G)$ will be evaluated as the trace of $G$ considered as an operator in ${\rm Hom}(\overline V_r, \overline V_r)$. A basis of $\overline V_r$ is given by the core circle $e_j$ of a solid torus bounded by ${\mathbf T}$, labeled with an {\em even} integer $0\leq j\leq r-2$. The result of $G$ applied to a basis element $e_j$ may be computed by pushing $G$ (considered as a quantum spin network with edges labeled by $2$) into the solid torus and re-expressing the result as a linear combination of $\{e_j\}$ using the recoupling theory \cite{KL}. In the examples below the trace will be computed using the expansion ${\Phi}(G)$ of $G$ in terms of multi-curves; the multi-curves in question will act diagonally on $\overline V_r$ with respect to the basis $\{ e_j, 0\leq j\leq r-2, j \; {\rm even} \}$. We emphasize that while multi-curves (consisting of ``spin $1/2$ loops'', or in the usual TQFT terminology loops labelled $1$) are elements of ${\rm SU}(2)$ and not ${\rm SO}(3)$ theory, the configurations of multi-loops considered below preserve the subspave $\overline V_5$, and they provide a convenient evaluation method.

We give several sample calculations of the trace ${\rm tr}^{}_5$, used in the proof below. First consider $k$ trivial loops (labelled $1$) on the torus. 
The usual $d$-isotopy relation states that removing a loop gives a factor $d={\phi}$, so ${\rm tr}^{}_5$ equals ${\phi}^k$ times the trace of the empty diagram. Since the dimension of the space (spanned by the core of the solid torus with labels $0$ and $2$) is $2$, the result is $2 {\phi}^k$. The trace evaluation of the {\em graph} consisting of $k$ trivial loops (figure \ref{trace11 fig}(a)) equals $2(d^2-1)^k$, which (precisely at this root of unity!) also equals $2 {\phi}^k$.

Next consider non-trivial (spin $1/2$) loops on the torus in (\ref{trace0 eq}). Using (\ref{LinkingLoop fig}), the action on $\overline V^{}_{5}$ is seen to be diagonal:
\begin{equation}  \label{trace0 eq}
\vcenter{\hbox{\includegraphics[height=2.3cm]{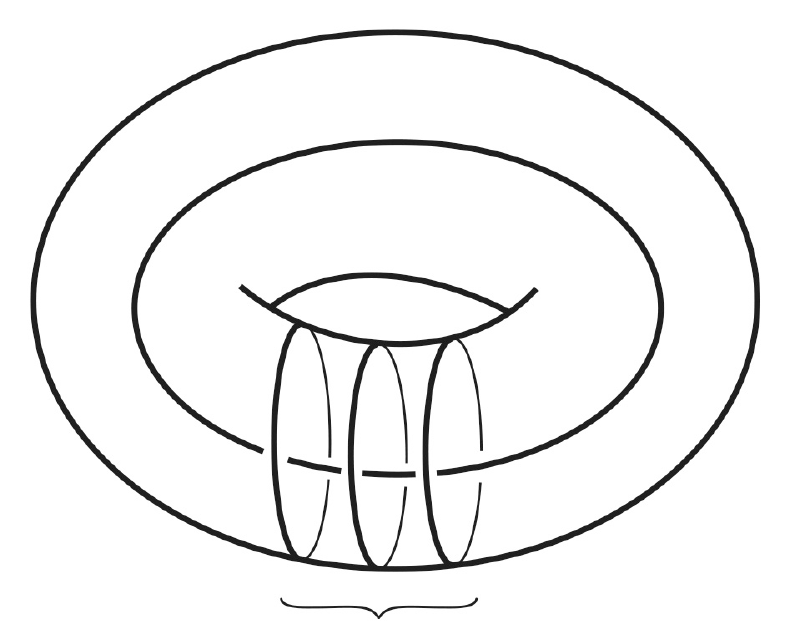}\hspace{4.2cm}  \includegraphics[height=2.3cm]{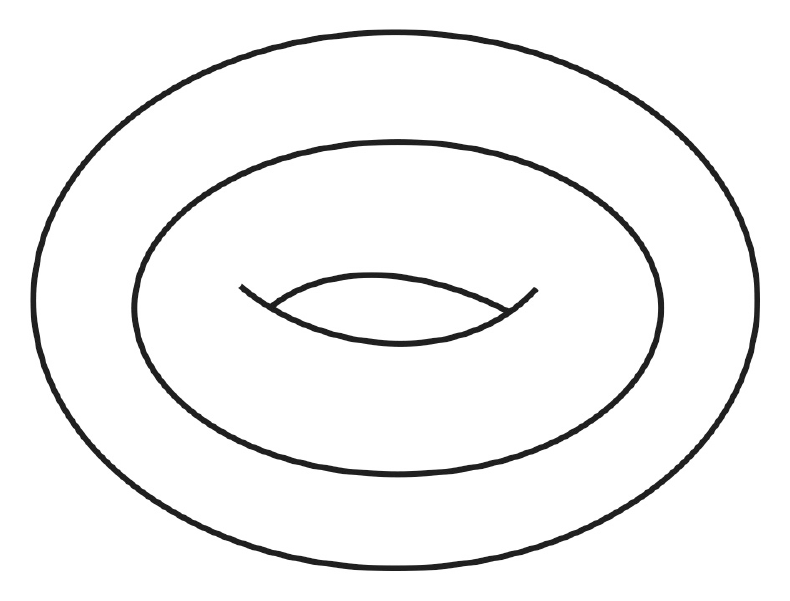}
\put(-255,-14){$k$}
{\scriptsize 
\put(-275,48){$j$}
\put(-68,48){$j$}
}
{\large 
\put(-202,30){$=\, \left[ \frac{{\rm sin}(2(j+1){\pi}/5)}{{\rm sin}((j+1){\pi}/5)} \right]^k$}
} 
}}
,
\end{equation}
$j=0,2$. Therefore the trace of $k$ non-trivial loops with label $1$ on the torus equals 
\begin{equation}  \label{trace1 eq}
\vcenter{\hbox{\includegraphics[height=2.3cm, trim=15cm 0 0 0]{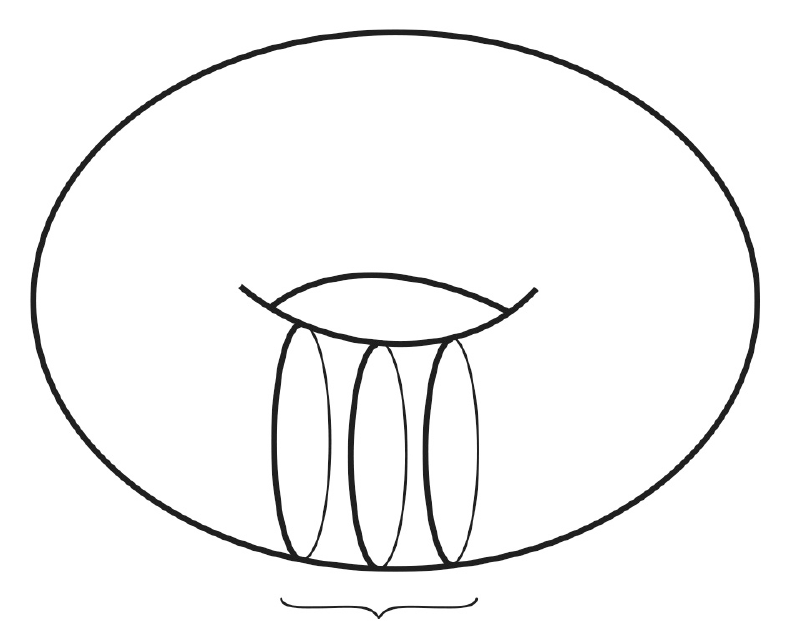}
\put(-47,-14){$k$}
{\large 
\put(-103,27){${\rm tr}^{}_{5}$}
\put(3,30){$=\, \left[ \frac{{\rm sin}(2{\pi}/5)}{{\rm sin}({\pi}/5)} \right]^k+\left[ \frac{{\rm sin}(6{\pi}/5)}{{\rm sin}(3{\pi}/5)} \right]^k$}}
\put(145,30){$=\, {\phi}^k+(-{\phi}^{-1})^k.$}

}}
\end{equation}

The analogous calculation for the {\em graph} consisting of $k$ non-trivial loops on the torus (or ``spin $1$ loops'') gives
\begin{equation}  \label{trace2 eq}
\vcenter{\hbox{\includegraphics[height=2.3cm, trim=15cm 0 0 0]{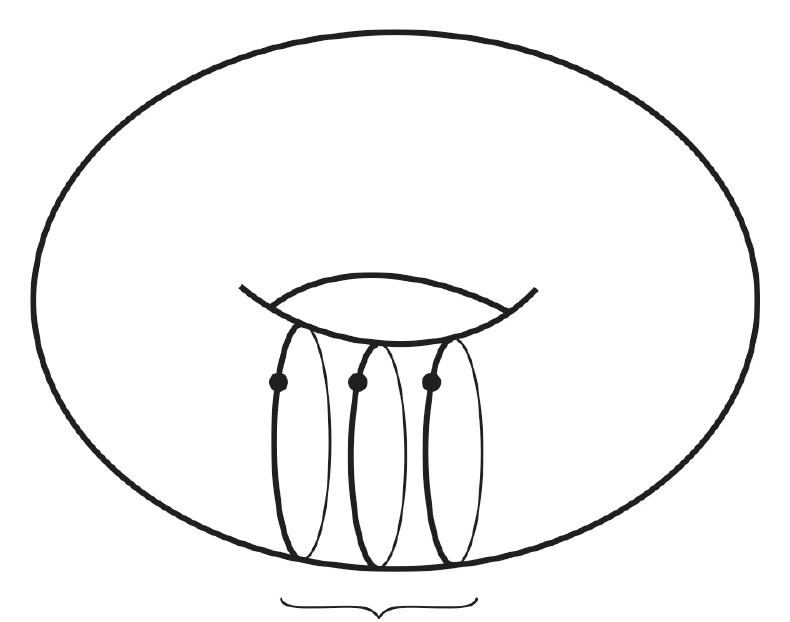}
\put(-47,-14){$k$}
{\large 
\put(-103,27){${\rm tr}^{}_{5}$}
\put(3,30){$=\, \left[ \frac{{\rm sin}(3{\pi}/5)}{{\rm sin}({\pi}/5)} \right]^k+\left[ \frac{{\rm sin}(9{\pi}/5)}{{\rm sin}(3{\pi}/5)} \right]^k$}}
\put(145,30){$={\phi}^k+(-{\phi}^{-1})^k.$}

}}
\end{equation}
The answer is again the same as for spin $1/2$ loops precisely at the $5$th root of unity. Note that the ${\rm SO}(3)$ trace is invariant under modular transformations of the torus, so (\ref{trace2 eq}) gives the trace of {\em any} $k$ non-trivial, spin $1$ loops on the torus. The situation is a bit more subtle with spin $1/2$ loops: the calculations in (\ref{trace0 eq}), (\ref{trace1 eq}) work specifically for non-trivial loops which bound disks intersecting the core of the solid torus once. A single spin $1/2$ loop which wraps in some other way around the torus and acts as a curve operator $V_5\longrightarrow V_5$, does not have to preserve the subspace $\overline V_5$. Nevertheless, an {\em even} number of non-trivial curves preserve $\overline V_5$, and moreover the evaluation (\ref{trace1 eq}) for $k$ even is in fact modular invariant: using (\ref{phi}), a pair of parallel spin $1/2$ loops may be expressed as a spin $1$ loop plus a scalar multiple of a trivial loop. This property will be used in the following section to evaluate the trace of graphs on the torus in terms of surround loops.

{\bf Remark.} There are two equivalent ways of computing the trace in (\ref{trace2 eq}): one using the formula  (\ref{LinkingLoop fig}) directly, or alternatively the $2$nd JW projectors can be expanded into linear combinations of spin $1/2$ loops, reducing the calculation to (\ref{trace1 eq}).

\begin{figure}[ht]
  \includegraphics[height=5.3cm, trim=-.9cm 0 0 0]{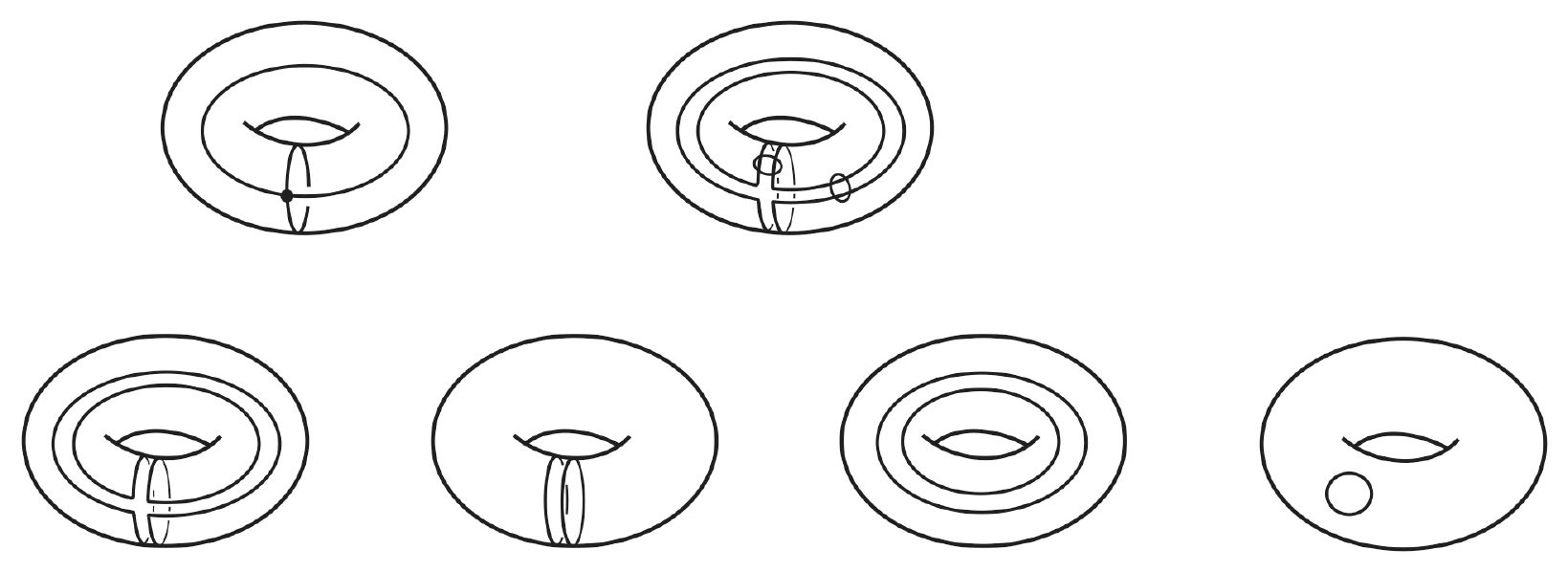}
\put(-390,114){${\rm tr}^{}_{5}$}
\put(-285,114){$= \;  d\, {\rm tr}^{}_{5}$}
\put(-155, 114){$=$}
\put(-435,31){$\; d{\rm tr}^{}_{5}$} 
\put(-328,31){$-\, {\rm tr}^{}_{5}$}
\put(-222,31){$-\, {\rm tr}^{}_{5}$}
\put(-116,31){$+\frac{1}{d}{\rm tr}^{}_{5}$}
\put(-320,-30){$=\, 2{\phi}^2-2({\phi}^2+(\frac{1}{\phi})^2)+2\, =\, 2-\frac{2}{{\phi}^2}$}
\caption{The factor $d$ in the top line comes from the normalization of the map $\Phi$ (figure \ref{phi}). At the $5$th root of unity $d={\phi}$. Both first and last terms in the bottom line have a trivial loop, evaluating to $d$, and the factor $2$ corresponds to the trace of the empty diagram. The two middle terms are evaluated according to (\ref{trace1 eq}).}
\label{trace3 fig}
\end{figure}

Finally, the trace of the graph in figure \ref{trace3 fig} is obtained by expanding both second projectors, and applying  (\ref{trace1 eq}).

\subsection{Trace and graph evaluations}
\begin{defi}   \label{poly eval def}
Given a graph $G\subset {\mathbb T}$, consider
\begin{equation} \label{P5 equ}
R^{}_5(G):=P^{}_G({\phi}^2, 1, {\phi}^{-2})+P^{}_G({\phi}^2, {\phi}^{-4}, {\phi}^{-2})
\end{equation}
\end{defi}
The two summands in the definition of $R^{}_5(G)$ are given by evaluations of the polynomial $P^{}_G$ in (\ref{poly definition}). In both cases, $Y={\phi}^2$ and $A={\phi}^{-2}$. Note that the first summand corresponds to $YW={\phi}^2$, and the second one to $YW={\phi}^{-2}$.

\begin{theoremsection} \label{equality at B5}
{\sl
Given any graph $G\subset {\mathbb T}$, the ${\rm SO}(3)$ TQFT trace evaluation ${\rm tr}^{}_{5}(G)$ equals $R^{}_5(G)$.
}
\end{theoremsection}

This result is the TQFT version on the torus of the loop evaluation of the flow polynomial of planar graphs at $Q={\phi}+1$, corresponding to $q=e^{2{\pi}i/5}$ (see (\ref{notation1})). 

{\em Proof of Theorem} \ref{equality at B5}.
We begin the proof by comparing calculations of ${\rm tr}^{}_{5}(G)$ and $R^{}_5(G)$ for the graphs in figure \ref{trace11 fig}, using results of Sections \ref{graph poly section}, \ref{trace subsection}.

(a) $G$ consists of $k$ trivial loops, figure \ref{trace11 fig}a.
$${\rm tr}^{}_{5}(G)=2(d^2-1)^k= 2{\phi}^k. \; \; R^{}_5(G)=(Y-1)^k+(Y-1)^k.$$
The factor $2$ in the expression for ${\rm tr}^{}_{5}(G)$ comes from the dimension $2$ of the vector space spanned by the even labels $0, 2$. Since $Y={\phi}^2$, the two expressions coincide.

(b) $k$ non-trivial loops, figure \ref{trace11 fig}b.  
According to (\ref{trace2 eq}), ${\rm tr}^{}_{5}(G)={\phi}^k+(-{\phi}^{-1})^k.$  By (\ref{trace12 fig}), $R^{}_5(G)=(YW-1)^k|_{YW={\phi}^2}+(YW-1)^k|_{YW={\phi}^{-2}}.$
Individual terms match:
$({\phi}^2-1)^k+({\phi}^{-2}-1)^k={\phi}^k+(-{\phi}^{-1})^k.$

(c) The graph in figure \ref{trace11 fig}c. The TQFT calculation in figure \ref{trace3 fig} gives 
 ${\rm tr}^{}_{5}(G)=2{\phi}^2-2({\phi}^2+(\frac{1}{\phi})^2)+2.$
By  (\ref{poly example}), $P^{}_G(Y, W, A)=AY^2- 2YW+1$. The two evaluations of $P^{}_G$, contributing to $R^{}_5(G)$, give ${\phi}^2-2{\phi}^2+1$ and ${\phi}^2-2(\frac{1}{\phi})^2+1$, adding up to the expression for ${\rm tr}^{}_{5}(G).$

The proof of Theorem \ref{equality at B5} for an arbitrary graph $G\subset {\mathbb T}$ is obtained by expanding the second JW projectors for all edges. The resulting summands for the TQFT trace are in $1$-$1$ correspondence with spanning subgraph $H\subset G$. To be precise, given a spanning subgraph $H$, in this correspondence the first term in the expansion (\ref{phi}) of the $2$nd projector is taken for each edge $e$ in $H$, and the second term is taken for each edge $e$ in $G\smallsetminus H$. The resulting loop configuration, called the {\em surround loops}, is the boundary of a regular neighborhood of $H$ on the surface. Each individual term in the trace evaluation equals
the sum of two entries, corresponding to the two labels $0, 2$, and we show next that  these entries precisely match the corresponding terms in the expansions 
 $P^{}_G({\phi}^2, 1, {\phi}^{-2})$, $P^{}_G({\phi}^2, {\phi}^{-4}, {\phi}^{-2})$.

For each spanning subgraph $H$ there are three cases, analogous to the examples (a)-(c) above:

(${\mathcal A}$) $H$ is homologically trivial on the torus: $H_1(H)\longrightarrow H_1({\mathbb T})$ is the zero map. The exponents of the variables $W$ and $A$ in (\ref{poly definition}) are zero in this case. The proof that loop evaluation in TQFT equals the summand in the definition of the graph polynomial (\ref{poly definition}) is thus identical to the planar case \cite[Lemma 2.5]{FK}. Both quantities in the statement of the theorem have a factor $2$: for ${\rm tr}^{}_{5}(G)$ this is because ${\rm dim}(\overline V^{}_5)=2$; for $R_{5}(G)$ the factor is the result of adding two identical summands in (\ref{P5 equ}).

(${\mathcal B}$) The image of $H_1(H)\longrightarrow H_1({\mathbb T})$ has rank $1$. In this case ${\rm s}(H)=0$. Recall that ${\bar {\rm c}(H)}$ denotes the number of connected components $H^{(i)}$ of $H$ such that $H_1(H^{(i)})\longrightarrow H_1({\mathbb T})$ of rank $1$ for each $i$. The term in the expansion of $R^{}_5(G)$ corresponding to $H$ is 
$$(-1)^{E(G)-E(H)} \,  {\phi}^{2{\rm n}(H)}\,  [1^{\bar {\rm c}(H)}+{\phi}^{-4{\bar {\rm c}(H)}}]\, = \, (-1)^{E(G)-E(H)} \,  {\phi}^{2({\rm n}(H)-\bar{\rm c}(H))}\,  [{\phi}^{2\bar {\rm c}(H)}+{\phi}^{-2\bar {\rm c}(H)}] .$$

Each $H^{(i)}$ has two surround loops which are non-trivial on the torus. 
In the calculation of ${\rm tr}^{}_5(G)$, these $\bar {\rm c}(H)$ non-trivial loops give a factor ${\phi}^{2\bar {\rm c}(H)}+{\phi}^{-2\bar {\rm c}(H)}$, matching the factor in square brackets in the calculation of $R^{}_5(G)$ above. 

The last step is to check that the remaining factor $(-1)^{E(G)-E(H)} \,  {\phi}^{2({\rm n}(H)-\bar{\rm c}(H))}$ above corresponds to the normalization and the trivial surround loops in the trace evaluation. As explained after (\ref{phi}), the normalization factor in the definition of $\Phi$ is $d^{E(G)-V(G)}$. Moreover, each edge in $G\smallsetminus H$ gives rise to an additional factor $-d^{-1}$ coming from the second term of the JW projector. This gives the desired sign $(-1)^{E(G)-E(H)}$. Thus the overall normalization factor corresponding to $H$ is $d^{E(H)-V(H)}$, where $d={\phi}$. In addition, each trivial surround loop of $H$ gives a factor $\phi$ in the trace evaluation. An Euler characteristic count gives the equality $$E(H)-V(H)+{\rm number\; of\; trivial\; loops}= 2({\rm n}(H)-\bar{\rm c}(H)),$$ concluding the proof in case ${\mathcal B}$.

(${\mathcal C}$) $H_1(H)\longrightarrow H_1({\mathbb T})$ is surjective, so ${\rm s}(H)=1$ and $\bar {\rm c}(H)=0$. This case is similar to (${\mathcal A}$) since all surround loops are trivial on the torus. Because of the homological assumption, there are two fewer surround loops than in the planar case, expected from the nullity ${\rm n}(H)$. In the TQFT evaluation this undercount gives a factor $d^{-2}$. This factor precisely matches the factor $A^{s(H)}=A={\phi}^{-2}$ in (\ref{poly definition}). \qed

Recall that the vanishing of the Jones-Wenzl projector is built into the definition of the TQFT vector space at the corresponding root of unity, so the $4$-th JW projector gives a local relation in $\overline V^{}_5$. 

\begin{cor} \label{corollary}
{\sl
The graph evaluation $R^{}_5(G)$ satisfies the local relation (\ref{fig:graphp4}), corresponding to the $4$th Jones-Wenzl projector. More precisely, given three graphs $G^{}_X, G^{}_I, G^{}_E$ on the torus, locally related as shown in figure (\ref{fig:graphp4}), $${\phi}R^{}_5(G^{}_X)=R^{}_5(G^{}_I)+R^{}_5(G^{}_E).$$
}
\end{cor}
More generally, given $G\subset {\mathbb T}$ and  odd $r$, consider 
\begin{equation} \label{Pr equ}
R^{}_r(G):=\sum_{j=0, \, j\, {\rm even}}^{r-2}  P^{}_G(d^2, W_{j,r} , d^{-2}),
\end{equation}
where $d=q^{1/2}+q^{-1/2}$ is the TQFT loop value corresponding to the root of unity $q=e^{2{\pi} i/r}$, and $W_{j,r}$ is defined by 
$$YW_{j,r}-1= \frac{{\rm sin}(2(j+1){\pi}/r)}{{\rm sin}((j+1){\pi}/r)}.$$

The proof of the following result is directly analogous to that  of Theorem \ref{equality at B5}, with TQFT sectors precisely corresponding to the summands in (\ref{Pr equ}):

\begin{theoremsection} \label{evaluations are equal} \sl
The ${\rm SO}(3)$ TQFT trace evaluation of $G$ at $q=e^{2{\pi}i/r}$  equals $R^{}_r(G)$. 
\end{theoremsection}

{\bf Remark.} A generalization of the polynomial $P$ for links $L$ in ${\mathbb T}\times [0,1]$ (along the lines of \cite[Section 6]{Kr}) gives a similar expression for the ${\rm SU}(2)$ trace of $L$. A polynomial $P_L$ for more general links in a surface $\Sigma$ times the circle was formulated in \cite{MS}. It is an interesting question whether the polynomial of \cite{MS} can be defined for ribbon {\em graphs} in ${\Sigma}\times S^1$, and whether our results extend to this setting.

\section{Golden identity for graphs on the torus}   \label{golden torus section}
Using TQFT methods developed above, in this section we formulate and prove an extension of the Tutte golden identity for graphs on the torus.

\subsection{Proof of the Tutte golden identity (\ref{golden identity eq})  for planar graphs} \label{planar subsection}

We start by summarizing the ideas underlying the proof in the planar case. Following \cite{FK}, we prove here the golden identity for the flow polynomial of cubic planar graphs $G$: $F^{}_G({\phi}+2) \, = \,  {\phi}^E\, (F^{}_G({\phi}+1))^2.$
The version of the contraction-deletion rule for cubic graphs reads
\begin{equation}  \label{trivalent figure}
\vcenter{\hbox{
\includegraphics[height=1.85cm]{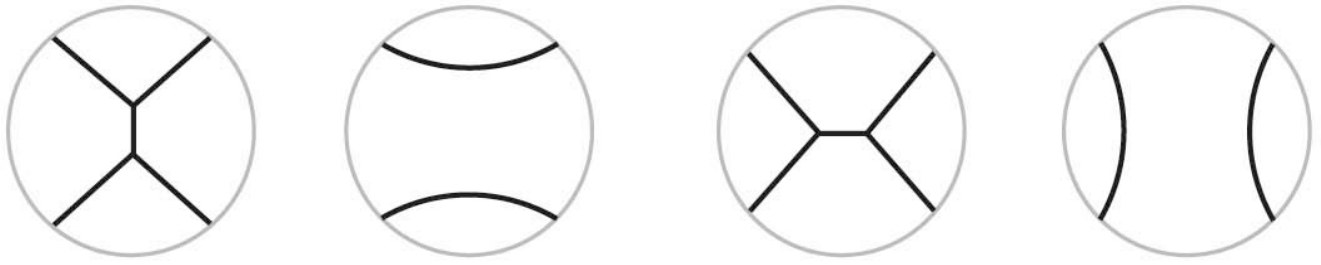}
    \put(-209,25){$+$}
    \put(-138,25){$=$}
    \put(-66,25){$+$}
 }}
\end{equation}
Using induction on the number of edges, it suffices to show that if three of the graphs in (\ref{trivalent figure}) for $F^{}_G({\phi}+2)$ satisfy the golden identity, then the fourth one does as well.
Given a cubic graph $G$, consider ${\phi}^E\, (F^{}_G({\phi}+1))^2$. It is convenient to formally depict, as in  (\ref{fig:tutte}), two identical copies of $G$, each one evaluated at ${\phi}+1$, with an overall factor ${\phi}^E= {\phi}^{3V/2}$. 
\begin{equation}  \label{fig:tutte}
\vcenter{\hbox{
\includegraphics[height=1.7cm]{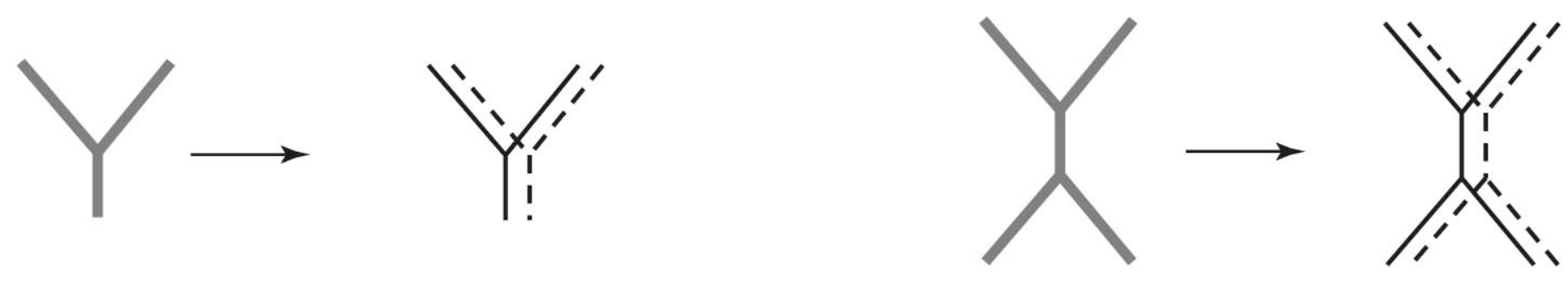}
{\scriptsize
    \put(-232,25){${\Psi}$}
    \put(-212,20){${\phi}^{3/2}\cdot$}
    \put(-62,26){${\Psi}$}
    \put(-40,20){${\phi}^3\cdot$}}
}}
\end{equation}
Here $E$ and $V$ denote the number of edges and vertices, respectively, of $G$. This doubling of lines is {\em not} that in the map $\Phi$ using the projector $p_2$, but instead a map ${\Psi}(G)=G\times G$.
Note that for $G=\,$circle,  $F_{\rm circle}({\phi}+2)={\phi}+1$. The corresponding value
for $\Psi$(circle) is ${\phi}^2={\phi}+1$, indeed the same.

\begin{figure}[ht]
\includegraphics[height=1.6cm]{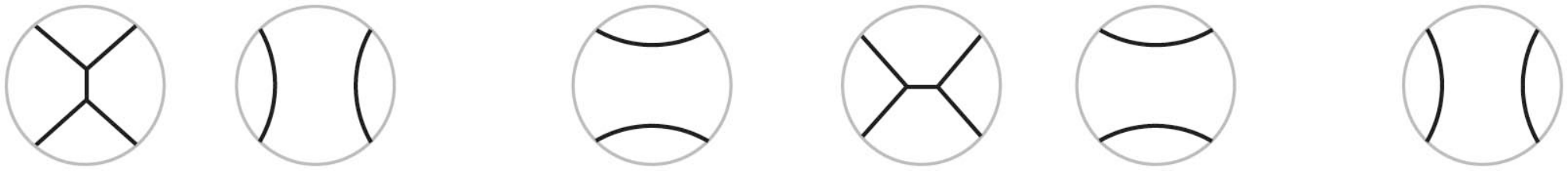}
{\small
    \put(-427,20){$\phi$}
    \put(-370,20){$=$}
    \put(-312,20){$+\,(1-{\phi})$}
    \put(-205,20){$\phi$}
    \put(-144,20){$=$}
    \put(-87,20){$+\,(1-{\phi})$}}
    \put(-223,7){$,$}
\caption{Local relations for the flow polynomial at $Q={\phi}+1$, equivalent to Tutte's relation (\ref{fig:graphp4}).}
\label{fig:JW}
\end{figure}

The strategy is to check that the evaluation ${\phi}^E\, (F^{}_G({\phi}+1))^2$ satisfies the relation (\ref{trivalent figure}) as a consequence of the {\em additional} local relation at $Q={\phi}+1$.
This additional relation is the graph version (\ref{fig:graphp4}) of the $4$th Jones-Wenzl projector.
Using the contraction-deletion rule, one checks  that (\ref{fig:graphp4}) is
equivalent to each of the  two relations shown in figure \ref{fig:JW}.

\begin{figure}[ht]
\includegraphics[height=1.55cm]{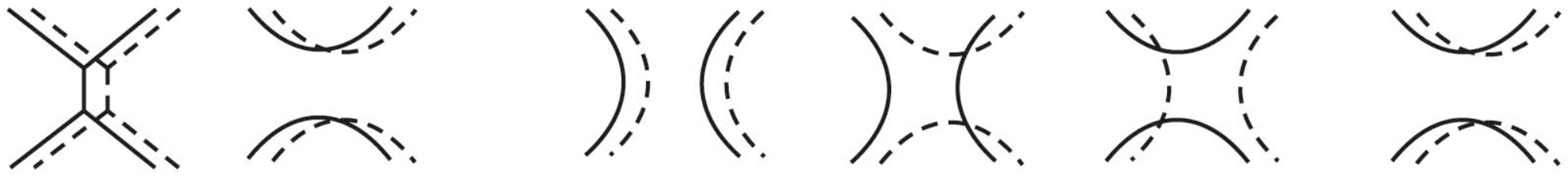}
    \put(-412,19){${\phi}^3\,\cdot$}
    \put(-348,19){$+$}
    \put(-283,19){$=\;{\phi}\,\cdot$}
    \put(-198,19){$-$}
    \put(-131,19){$-$}
    \put(-69,19){$+\;{\phi}\,\cdot$}
\caption{}
\label{tutte3}
\end{figure}

Consider the image of (\ref{trivalent figure}) under $\Psi$:
\begin{equation}  \label{tutte2}
\vcenter{\hbox{
\includegraphics[height=1.65cm]{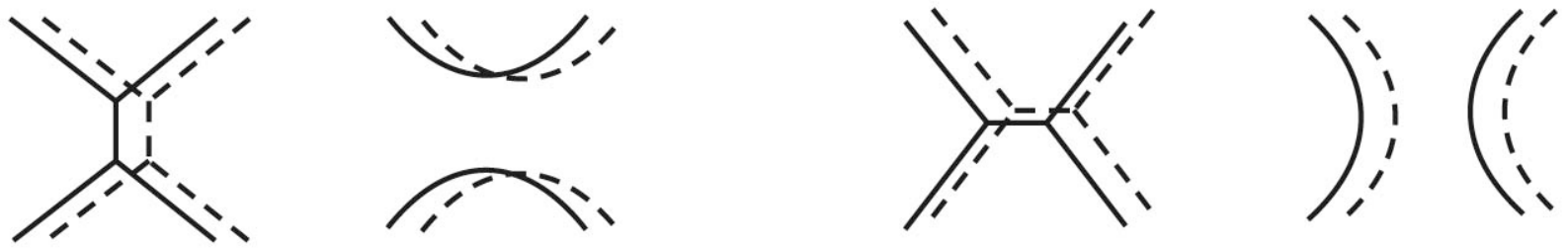}
    \put(-310,21){${\phi}^3\, \cdot$}
    \put(-241,21){$+$}
    \put(-166,21){$= \;{\phi}^3\, \cdot$}
    \put(-70,21){$+$}
    }}
\end{equation}
Applying the relation on the left in figure \ref{fig:JW} to both copies of the graph on the left in figure (\ref{tutte2}) yields figure (\ref{tutte3}).
The resulting expression on the right is invariant under 90 degree rotation, so must also be equal to the graph on the right of (\ref{tutte2}). Thus (\ref{tutte2}) holds, showing that the evaluation ${\phi}^E\, (F^{}_G({\phi}+1))^2$ satisfies the contraction-deletion relation (\ref{trivalent figure}). This concludes the proof of the golden identity for the flow polynomial of planar cubic graphs. \qed

\subsection{Extension to graphs on the torus}
The expression (\ref{Pr equ}) is defined for odd $r$. We generalize it with the following sum of evaluations of the graph polynomial $P^{}_G$ in (\ref{poly definition}):
\begin{equation} \label{P10 equ}
R^{}_{10}(G):= P^{}_G(Y, W_1, A)+P^{}_G(Y, W_2, A)+2P^{}_G(Y, W_3, A),
\end{equation}
where $Y={\phi}+2$, $A =(2({\phi}+2))^{-1}$, and the values $W_j$, $j=1,2,3$, are defined by 
$$YW_1={\phi}+2, \;  YW_2=1+{\phi}^{-2}, \; YW_3=0.$$
The choice of these values may be thought of as a choice of particular sectors of the ${\rm SO}(3)$ TQFT vector space of the torus at $q=e^{2{\pi}i/10}$. 
(Our forthcoming work, relating these results to lattice models on the torus, will give further evidence for why this is a relevant invariant at this root of unity.) 
The value of $R^{}_{10}(G)$ for $G$ a trivial loop on the torus equals $4(Y-1)=4{\phi}^2$.
Using (\ref{trace12 fig}), the value of $R^{}_{10}(G)$ for the graph consisting of a single non-trivial loop is seen to be ${\phi}^2+{\phi}^{-2}-2$. 

We are in a position to state the main result of this section.

\begin{theoremsection} \label{golden torus theorem} \sl
Let $G\subset {\mathbb T}$ be a cubic graph. Then 
\begin{equation}   \label{golden torus eq}
R^{}_{10}(G)\, =\, {\phi}^E  R^{}_{5}(G)^2,
\end{equation}
where $E$ is the number of edges of $G$.
\end{theoremsection}

The version of the contraction-deletion relation for cubic graphs is shown in   (\ref{trivalent figure}). As usual, we allow cubic graphs with disjoint loops; such loops do not count towards $E$.
The proof in the planar case (see section \ref{planar subsection}) showed that if three of the graphs in (\ref{trivalent figure}) satisfy the golden identity, then the fourth one satisfies it as well. This fact holds for the identity \eqref{golden torus eq} for graphs on the torus as well. Indeed, (\ref{trivalent figure}) holds for $R^{}_{10}(G)$ since it is defined as the sum (\ref{P10 equ}) of polynomials satisfying the contraction-deletion rule. And (\ref{trivalent figure}) holds for ${\phi}^E  R^{}_{5}(G)^2$ for the same reason as in section \ref{planar subsection}, since by Corollary \ref{corollary} $R^{}_{5}(G)$ obeys the local relation corresponding to the $4th$ JW projector.

For cubic graphs $G$ which are homologically trivial on the torus  
($r(G\subset {\mathbb T}) =0$ in the notation of section \ref{graph poly section}) the proof of (\ref{golden torus eq}) follows from the planar case in section \ref{planar subsection}
since the polynomial $P^{}_G$ in (\ref{poly definition}) equals the flow poynomial $F^{}_G$. (For example, calculations in section \ref{graph poly section} show that in the special case of the graph consisting of $k$ trivial loops, $R^{}_5(G)=2{\phi}^k$, and $R^{}_{10}(G)= 4{\phi}^{2k}=R^{}_5(G)^2$.)

Next consider the case of cubic graphs $G$ of rank $r(G\subset {\mathbb T}) =2$. Consider a minimal cubic graph $G$ of this type,  shown in (\ref{fig:TorusGraph}),  where the square with opposite sides identified is the usual representation of the torus. A direct calculation using (\ref{poly definition}), or alternatively using the contraction-deletion rule to reduce this to calculations above, shows 
 $R^{}_5(G)=R^{}_{10}(G)={\phi}^{-3}$, proving (\ref{golden torus eq}) in this case.
\begin{equation} \label{fig:TorusGraph}
\vcenter{\hbox{        
\includegraphics[height=2.7cm]{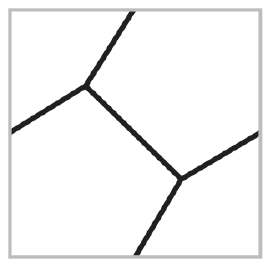}
}}
\end{equation}

The proof in general for rank $2$ cubic graphs is by induction on the number of edges.
It is proved in \cite{D} that any two triangulations of the torus with the same number of vertices are related by diagonal flips, up to equivalence given by diffeomorphisms.
(This was extended in \cite{Negami} to pseudo-triangulations of surfaces of any genus, where an embedding $\Gamma\subset {\mathbb T}$ is a {\em pseudo-triangulation} if each face is a three-edged $2$-cell, possibly with multiple edges and loops.) Formulated in terms of dual cubic graphs,  two cellular embeddings of cubic graphs (or in other words rank $2$ graphs) on the torus are related by the $I-H$ move:
\begin{equation} \label{fig:flip}
\vcenter{\hbox{
\includegraphics[height=2cm]{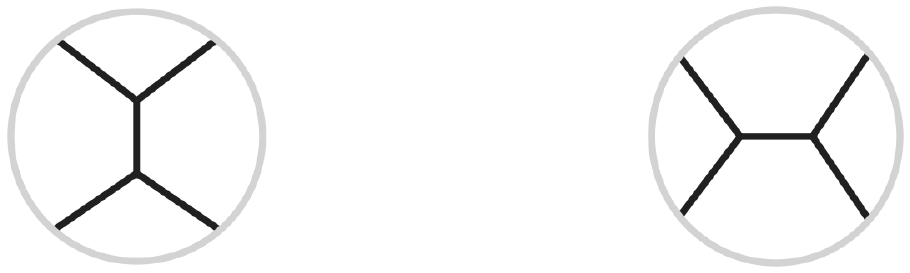}
 }}
\end{equation}
The relation (\ref{trivalent figure}) accomplishes the $I-H$ move, while also introducing graphs with fewer edges. The theorem has been checked for a minimal cubic graph in  (\ref{fig:TorusGraph}), and the inductive step is achieved by a local modification (\ref{fig:local}), which increases the number of edges and preserves (\ref{golden torus eq}). 
\begin{equation} \label{fig:local}
\vcenter{\hbox{
\includegraphics[height=1.9cm]{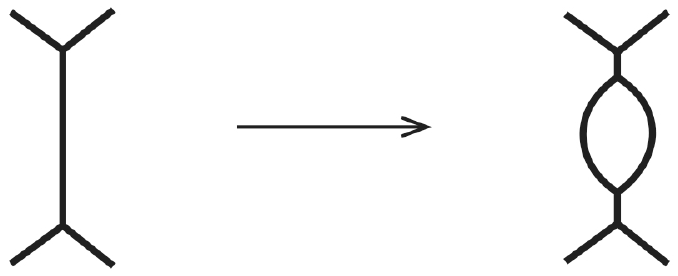}
}}
\end{equation} 

Finally consider  cubic graphs of rank $1$ (that is, $r(G\subset {\mathbb T}) =1$), or equivalently graphs on the cylinder.
The fact that triangulations of the sphere with the same number of vertices are related by diagonal flips dates back to \cite{Wagner}. 
Using (\ref{trivalent figure}) and (\ref{fig:local}) as above, the proof in the rank $1$ case therefore follows from the calculation for $G$ consisting of  $k$ non-trivial loops on the torus. 
In this case, using (\ref{trace12 fig}) one has $R^{}_5(G)={\phi}^k+(-{\phi}^{-1})^{k}$, and $R^{}_{10}(G)= {\phi}^{2k}+{\phi}^{-2k}+2(-1)^k=R^{}_5(G)^2$. 
\qed

{\bf Remark.} Theorem \ref{golden torus theorem} is an extension of the golden identity for the flow polynomial of planar cubic graphs. The focus of this paper is on the ``topological flow polynomial'' $P^{}_G$, since it is related to the TQFT trace evaluation, as stated in Theorem \ref{evaluations are equal}.  It is worth noting that an analogue of the invariants $R^{}_5(G), R^{}_{10}(G)$ may be defined using the polynomial  $C^{}_G$ in place of $P^{}_G$. Using the duality relation (\ref{duality eq}), the identity (\ref{golden torus eq}) then gives rise to an extension of the original Tutte golden identity (\ref{golden identity eq}) for the ``topological chromatic polynomial'' $C^{}_G$ of triangulations of the torus.

{\bf Acknowledgments.} VK would like to thank Calvin McPhail-Snyder and Sittipong Thamrongpairoj for discussions about polynomial invariants of graphs on surfaces.
 VK was supported in part by NSF grant DMS-1612159; he also would like to thank All Souls College and the Department of Physics at the University of Oxford for hospitality and support.


\begin{thebibliography}{10}

\bibitem{AK} I. Agol and V. Krushkal,
{\em Tutte relations, TQFT, and planarity of cubic graphs}, 
Illinois J. Math. 60 (2016), no. 1, 273-288. 

\bibitem{AK2} I. Agol and V. Krushkal,
{\em Structure of the flow and Yamada polynomials of cubic graphs},
 arXiv:1801.00502.

\bibitem{AMF}
  D.~Aasen, R.~S.~K.~Mong and P.~Fendley,
  {\em Topological Defects on the Lattice I: The Ising model},
  J.\ Phys.\ A {\bf 49} (2016) no.35,  354001

\bibitem{BHMV} C. Blanchet, N. Habegger, G. Masbaum and P. Vogel, {Topological quantum field theories
derived from the Kauffman bracket}, Topology 34 (1995), 883-927.

\bibitem{BR} B. Bollob\'{a}s and O.  Riordan, {\em A polynomial of graphs on surfaces}, Math. Ann. 323 (2002), no. 1, 81--96.

\bibitem{CR}  D. Cimasoni and N. Reshetikhin, {\em Dimers on surface graphs and spin structures. I}, Comm. Math. Phys. 275 (2007), 187--208.


\bibitem{DFKLS} O. Dasbach, D. Futer, E. Kalfagianni, X.-S. Lin, Xiao-Song and N. Stoltzfus, 
{\em The Jones polynomial and graphs on surfaces},
J. Combin. Theory Ser. B 98 (2008), 384--399.

\bibitem{D} A.K. Dewdney, 
{\em }Wagner's theorem for torus graphs,
Discrete Math. 4 (1973), 139--149.

\bibitem{DSZ} P. di Francesco, H. Saleur and J.-B. Zuber, {\em	
Relations Between The Coulomb Gas Picture And Conformal Invariance Of Two-dimensional Critical Models}, J. Stat. Phys. 49 (1987) 57



\bibitem{FK} P. Fendley and V.  Krushkal, {\em Tutte chromatic identities from the Temperley-Lieb algebra}, Geom. Topol. 13 (2009), 709--741.

\bibitem{FK2} P. Fendley and V.  Krushkal, {\em Link invariants, the chromatic polynomial and the Potts model}, Adv. Theor. Math. Phys.  14 (2010) 2, 507--540.

\bibitem{FK3} P. Fendley and V.  Krushkal, {\em Loop models and a topological Tutte polynomial for graphs on the torus}, in preparation.



\bibitem{FNWW} M. Freedman, C. Nayak,  K. Walker and Z. Wang, {\em On picture (2+1)-TQFTs}, Topology and physics, 19-106, Nankai Tracts Math., 12, World Sci. Publ., Hackensack, NJ, 2008. 

\bibitem{Jo} V.F.R.~Jones,
Subfactors and knots. 
CBMS Regional Conference Series in Mathematics, 80. Published for the Conference Board of the Mathematical Sciences, Washington, DC; American Mathematical Society, Providence, RI, 1991.

\bibitem{KL} L.H. Kauffman and S.L. Lins,
 Temperley-Lieb recoupling theory and invariants of 3-manifolds. 
Annals of Mathematics Studies, 134. Princeton University Press, Princeton, NJ, 1994.

\bibitem{KRTW} T. Krajewski, V. Rivasseau, A. Tanasă and Z. Wang, 
{\em Topological graph polynomials and quantum field theory. I. Heat kernel theories},
J. Noncommut. Geom. 4 (2010), 29--82.

\bibitem{Kr} V. Krushkal, {\em Graphs, links, and duality on surfaces}, Combin. Probab. Comput. 20 (2011), 267--287.

\bibitem{MS} J. March\'{e} and  R. Santharoubane, 
{\em Asymptotics of quantum representations of surface groups}, arXiv:1607.00664

\bibitem{MM}
C. McPhail-Snyder and K.A. Miller,  
{\em Planar diagrams for local invariants of graphs in surfaces}, 
arXiv:1805.00575.


\bibitem{MPS} S. Morrison, E. Peters and N. Snyder, {\em  Knot polynomial identities and quantum group coincidences}, Quantum Topol. 2 (2011), 101--156.

\bibitem{Negami} S. Negami, 
{\em Diagonal flips in pseudo-triangulations on closed surfaces}, 
Discrete Math. 240 (2001), 187--196.  



\bibitem{Pasquier} V. Pasquier, {\em  Lattice derivation of modular invariant partition function on the torus}, Journal of Physics A 20 L1229--L1237




\bibitem{T1} W.T. Tutte, {\em On chromatic polynomials and the golden ratio}, J. Combinatorial Theory 9 (1970), 289--296.


\bibitem{T2} W.T.~Tutte, {\em More about chromatic polynomials and the golden ratio}, Combinatorial Structures and their Applications, 439--453 (Proc. Calgary Internat. Conf., Calgary, Alta., 1969)


\bibitem{Wagner} K. Wagner, {\em Bemekungen zum Vierfarbenproblem}, J. Der Deut. Math. Ver. Abt. 1 46 (1936) 26--32.


\bibitem{We} H.Wenzl,
{\em On a sequence of projections},  C. R. Math. Rep. Acad. Sci. Canada  9  (1987), 5--9.





\end{thebibliography}
\end{document}